\documentclass{amsart}
\pdfoutput=1
\usepackage{amsmath,amscd,amssymb,amsthm,epsfig,psfrag,latexsym,subfigure,hyperref, mathrsfs, color}

\newtheorem{theorem}{Theorem}[section]

\newtheorem{conjecture}[theorem]{Conjecture}

\newcommand{\ZZ}{\mathbb{Z}}

\newcommand{\CC}{\mathbb{C}}

\newcommand{\co}{\colon\thinspace}

\def\vol{\mbox{\rm{Vol}}}

\def\min{\mbox{\rm{min}}}

\edef\t@mp{\catcode`\noexpand\#=\the\catcode`\#}%
    \catcode`\#=12
    \def\h@sh{#}%
    \t@mp

\edef\t@mp{\catcode`\noexpand\~=\the\catcode`\~}%
    \catcode`\~=12
    \def\tild@{~}%
    \t@mp

\begin{document}

\title{Problems In Groups, Geometry, and Three-Manifolds}
\author{Edited by Kelly Delp, Diane Hoffoss, and Jason Fox Manning}
\maketitle
\centerline{\large Dedicated to Daryl Cooper, on the occasion of his 60th birthday.}

\bigskip

In May 2015, a conference entitled ``Groups, Geometry, and $3$--manifolds'' was held at the University of California, Berkeley.\footnote{This conference was supported by NSF grant DMS-1503955.}  The organizers asked participants to suggest problems and open questions, related in some way to the subject of the conference.  These have been collected here, roughly divided by topic.  The name (or names) attached to each question is that of the proposer, though many of the questions have been asked before.

\setcounter{tocdepth}{1}
\tableofcontents

\section{Convex Projective Structures}
%%%%%%%%%%%%%%%%%%%%%%%
\subsection{} (Agol) If $M^n$ is a closed manifold with a strictly convex projective structure, is it cubulated?
A theorem of Benoist implies that $\pi_1(M)$ is hyperbolic \cite{Benoist01}. This question is open even in the
case of closed hyperbolic $n$-manifolds. 

\subsection{} (Choi) Suppose a hyperbolic $3$--manifold admits a CR-structure (not necessarily spherical).  Can the deformation theory of convex real projective structures on $M$ be understood in terms of the CR-structure?  Possibly this is easier for hyperbolic Coxeter $3$--orbifolds.

\subsection{} (Cooper) If $M$ is a closed hyperbolic 3-manifold is every projective structure on $M$ convex?

\subsection{} (Danciger) Let $N$ be the closed three-manifold obtained by gluing together two copies of figure eight knot complement along their torus boundaries by some homeomorphism. Does $N$ admit a convex projective structure?  (It follows from Ballas--Danciger--Lee \cite{BDL15} that there is a convex projective structure on the \emph{double} of the figure eight knot complement, i.e. if the homeomorphism is chosen to be the identity.)

\subsection{} (Danciger) More general version:
Let $N$ be a closed three-manifold whose JSJ decomposition contains only hyperbolic pieces. Does $N$ admit a convex projective structure?

%%%%%%end convex%%%%%%%%%%%%%
\section{Geometric Transitions}
%%%%%%%%%%%%%%%%%%%%%%%

\subsection{} (Leitner) Geometric transitions are continuous paths of geometries which abruptly change type in the limit (Cooper-Danciger-Wienhard \cite{CooperDW14}).  The most intuitive example is a sequence of spheres with increasing radius which limit to a plane.  It remains to understand all of the transitions between the eight Thurston geometries.  More generally, how can one tell when one geometry is a limit of another?  One can find a path of transitions to show one geometry limits to another, but it is in general much more difficult to show that one geometry can NOT limit to another.   What properties must be satisfied by a limiting geometry?

\subsection{} (Cooper) Does anyone know of a non-zero polynomial on the tensor product of vector spaces $U \otimes V \otimes W$ that is invariant under the action of $SL(U) \times SL(V) \times SL(W)$ when $\mbox{dim}(U)=\mbox{dim}(V)= 4$ and $\mbox{dim}(W)= 8$?  This is related to limits under conjugacy of the diagonal subgroup in $SL(8,\mathbb{R})$.

\subsection{} (Cooper) Suppose $\beta$ is  a non-degenerate bilinear form on a finite dimensional real vector space $V$ and $G = \mbox{Isom}(\beta) \subset GL(V)$.  Which subgroups $H \subset GL(V)$ are the Hausdorff limits of  sequences of conjugates of G?  This is known for non-degenerate symmetric and skew-symmetric forms.

% \subsection{} (Wienhard) $(\rho_{t}, \rho_{t+\epsilon}) : \Gamma \rightarrow O(p,q) \times O(q,p) \rightarrow O(p+q, q+p)$ (from Fanny Kassel's talk). Is it $P_1$ Anosov?\TODO{Add context to Wienhard questions or delete.}

% \subsection{} (Wienhard) If $\Gamma \rightarrow O(p,q) \times O(q,p) \rightarrow O(p+q, q+p)$ is $P_{p+1}$-Anosov, then $\Gamma \curvearrowright \mathbb{H}^{p,q}$ properly discontinuously \cite{GGKW15}.  Is there a converse to this?

%%%%%%%%%%%%%end transitions%%%%%%%%%%%%%%%

\section{Surfaces in 3-Manifolds}
%%%%%%%%%%%%%%%%%%%%%%%
\subsection{} (Cooper) Are hyperbolic 3-manifolds with large enough injectivity radius Haken \cite[Problem 3.58]{Kirby}? 

\subsection{}
(Agol) Given a 3-manifold $M$ whose fundamental group acts on a 
simplicial tree without global fixed points, consider the covering space
associated to the subgroup fixing an edge of the tree. This covering
space has a unique 2-dimensional homology class which separates
the two ends corresponding to the two ends of the tree minus the edge.
Is there a surface in this homology class with minimal Thurston norm that embeds in $M$?

\subsection{} (Agol) Does every closed hyperbolic 3-manifold contain a closed $\pi_1$-injective surface
with only double curves of intersection? It seems like this ought to be true 
for all but finitely many hyperbolic 3-manifolds with volume $< V$ for any $V$ by extending
results of \cite{cooperlong, TaoLi}.

\subsection{} (Agol) Does every closed hyperbolic 3-manifold contain a closed $\pi_1$-injective
surface which satisfies the 1-line property? This means that in the universal cover, any
pair of preimages of the surface intersect in a single line (and the intersection of stabilizers
is $\mathbb{Z}$). The surfaces constructed by Kahn and Markovic \cite{KM12} will likely not have
this property, since they tend to overlap on large subsurfaces. 

\subsection{} (Agol) Do cusped finite-volume hyperbolic 3-manifolds have closed quasi-fuchsian surfaces which
cubulate except for the cusps? Geometrically, for every pair of points in $\partial_\infty\mathbb{H}^3$,
some lift of the limit set of a quasi-fuchsian surface should separate this pair. This might be accessible
using constructions of Masters-Zhang \cite{MastersZhang08,MastersZhang09}  and Baker-Cooper \cite{BakerCooper14}. 
One would obtain a cocompact action on a CAT(0) cube complex in which the only point
stabilizers are parabolic subgroups \cite{BergeronWise12}. 

\subsection{} (Agol) Are finite-volume hyperbolic 3-manifolds virtually semi-fibered? 

\subsection{} (Agol) Which Kleinian groups admit closed quasi-fuchsian surface subgroups? The difficult case is
when the group has parabolic subgroups (when there are no parabolics, then this is true if and only if
the group is not virtually free). 

\subsection{} (Agol) Let $M$ be a $3$--manifold, $\phi\co \pi_1M\to \ZZ$ dual to the surface $(\Sigma,\partial\Sigma)\subset (M,\partial M)$.  Let $N$ be the $3$--manifold obtained by cutting $M$ along $\Sigma$.  For what $\alpha$ in $\mathrm{Hom}(\pi_1M,SL_2(\CC))$ is $N$ an $\alpha$--twisted homology product?

  A conjecture of Dunfield--Friedl--Jackson would say that $N$ is an $\alpha$--twisted homology product whenever $M$ is hyperbolic and $\alpha$ is discrete faithful.

\subsection{} (Futer, Manning) Hass's version of the Simple Loop Conjecture \cite[Problem 3.96]{Kirby} states:
\begin{conjecture}
  If $\phi\co\Sigma^2\to M^3$ is any two-sided immersion of a closed
surface of genus $\geq 1$, then either $\phi$ is $\pi_1$--injective or
\emph{compressible}. ($\phi$ is compressible if there is an essential
simple closed curve $\gamma$ on $\Sigma$ so that $\phi(\gamma)$ is
null-homotopic in $M$.)
\end{conjecture}
As pointed out by Hass, the conjecture is true if and only if it is true for irreducible $M$. The conjecture is known for $M$ a Seifert fibered space (Hass \cite{Hass87}), a graph manifold (Rubinstein-Wang \cite{RubinsteinWang98}), or a Solv-manifold (recent work of Drew Zemke \cite{zemke2015simple}).  It seems not to be known for any hyperbolic manifold, or for any $3$--manifold with a nontrivial JSJ with a hyperbolic piece.

Problem:  Prove it for some particular hyperbolic manifolds, or find a counterexample.

\subsection{} (Futer, Schleimer) Is there a practical algorithm to test whether a pair of 3-manifolds are homeomorphic?

\subsection{} (Walsh) The CAT(0) cubed dimension of a group $G$ is the
minimal dimension of a CAT(0) cubed space on which $G$ acts
geometrically. Is there a sequence of closed hyperbolic 3-manifold
groups, such that their CAT(0) cubed dimension tends to infinity? Is there some other sequence of CAT(0) groups such that the difference between the CAT(0) dimension and the CAT(0) cubed dimension tends to infinity?  Note, by Bridson \cite{Bridson01}, there are CAT(0) groups such that the difference between the geometric dimension and the CAT(0) dimension is arbitrarily large.

%%%%%%%%end surfaces%%%%%%%%%%

\section{Hyperbolic Groups and groups acting on spaces}

\subsection{} (Agol) Do closed hyperbolic 3-manifold groups have finite-index subgroups embedding in a word-hyperbolic
reflection group?

\subsection{} (Futer) Gromov asked whether 
every freely indecomposable hyperbolic group contains the fundamental group of a closed hyperbolic surface.  That the answer is yes is sometimes called the \emph{Surface Subgroup Conjecture} for hyperbolic groups.  Is this any easier for cubulated hyperbolic groups?

\subsection{} (Cooper) If $M$ is a closed 3-manifold when does $\pi_1 M$ act on the affine building for $SL(4,\mathbb{R})$ so that the quotient retracts to $M$?  An example is $M = \vol3$.

\subsection{} (Kassel, Mann) Let $\Gamma$ be a discrete group acting properly discontinuously by affine transformations on $\mathbb{R}^5$. Is $\Gamma$ virtually an extension of a free group by a solvable group?

\subsection{} (Kassel, Mann) Classify all properly discontinuous affine actions of a given closed surface group on $\mathbb{R}^6$.

\subsection{} (Kassel, Mann) What is the minimal $n$ such that a given a right-angled Coxeter group admits proper affine action on $\mathbb{R}^n$?

%%%%%end hyperbolic groups%%%%%%%%%

\section{Hyperbolic Geometry/Kleinian Groups}
%%%%%%%%%%%%%%%%%%%%%%%%%
\subsection{} (Agol) Characterize hyperbolic 3-manifolds with infinitely generated fundamental
group. In particular, is there a 3-manifold which is locally hyperbolic with 
no infinitely-divisible subgroup of the fundamental group (like $\mathbb{Q}$) but not hyperbolic?
By locally hyperbolic, we mean any cover with finitely generated fundamental
group admits a complete hyperbolic metric, and in particular is tame.

\subsection{} (Agol) Do there exist fibered hyperbolic 3-manifolds which are homology $S^2\times S^1$,
and have arbitrarily large injectivity radius? Are there hyperbolic homology spheres
of arbitrarily large injectivity radius? See \cite{BergeronVenkatesh10}, \cite{BrockDunfield15}. (cf. \cite{BostonEllenberg,CalegariDunfield} for the case of \emph{rational} homology spheres)

\subsection{} (Agol) Characterize the polytopes of the Thurston norms
of finite volume hyperbolic 3-manifolds. Thurston completed the rank 2 case, but
in higher rank, this is completely open \cite{Th3}. 

\subsection{} (Agol) Do hyperbolic 3-manifolds have a finite-sheeted cover homeomorphic to a CAT(0)-cube
complex? A good test case might be arithmetic 3-manifolds containing a geodesic surface. 

\subsection{} (Agol) Fix a constant $\mu$ less than the 3-dimensional Margulis constant. Consider the
hyperbolic 3-manifolds with volume $< V$, and drill out all closed geodesics of length $< \mu$. 
Of the resulting finite collection of 3-manifolds (throwing out repeats), let $s(V)$ be the fraction
of these cusped hyperbolic manifolds which are small (contain no closed incompressible non boundary parallel surface). What
is the limiting behavior of $s(V)$ as $V\to \infty$? How does it depend on $\mu$?

\subsection{} (Agol) If $M_1, M_2$ are cusped hyperbolic 3-manifolds, does there exist a cover $M_1'\to M_1$ and
a non-zero degree map $M_1'\to M_2$ taking cusps to cusps? 

\subsection{} (Agol) The renormalized volume of quasi-fuchsian groups was defined by 
Schlenker \cite{Schlenker13}. By Bers' simultaneous uniformization theorem, this gives a function $\rho: \mathcal{T}(S)\times \mathcal{T}(S)\to \mathbb{R}$,
where $\mathcal{T}(S)$ denotes the Teichumuller space of a surface $S$. Is $\rho$ a metric on $\mathcal{T}(S)$?

\subsection{} (Schleimer) Give a rigorous explanation for why SnapPy works so well in practice.

\subsection{} (Cooper)
Given $R>0$ and an integer $n\ge 4$ is there $\epsilon>0$ and
 a finite set of hyperbolic $n$-simplices with the following property. 
  Consider those  closed $n$-manifolds which can be constructed by gluing copies 
  of these simplices together by isometries along codimension-1 faces, such the resulting 
  hyperbolic cone $n$-manifold has cone angles around each codimension-2
simplex that are in the interval $(2\pi-\epsilon,2\pi+\epsilon).$ Then every such closed
 $n$-manifold admits a hyperbolic metric, and every closed hyperbolic $n$-manifold with injectivity radius everywhere bigger than $R$ is obtained in this way.  A set $P$ with these properties is called a {\em Thurston's Lego set  in $n$ dimensions}.
 These exist for $n\le 3$.

\subsection{} (Futer, Schleimer) The Triangulation Conjecture: Every cusped hyperbolic 3-manifold admits a geometric triangulation.

\subsection{} (Futer) Build a combinatorial model for hyperbolic 3-manifolds, with explicit bi-Lipchitz constants. 
%\TODO{Can you be more specific?  Do you mean finite volume?}

\subsection{} (Reid) Let $\Gamma$ be a Kleinian group of finite covolume, and let $\mathcal{C}(Γ\Gamma)$ be the set of isomorphism classes of the finite groups that are quotients (homomorphic images) of $\Gamma$.  Does $\mathcal{C}(\Gamma)$ determine $\Gamma$ up to isomorphism?  (see \cite{Reid13}, \cite{BridsonConderReid15})

\subsection{} (Reid) Let $\Gamma = F_r$, $r \ge 2$.  Does $\mathcal{C}(F_r)$ determine $F_r$ up to isomorphism?  

\subsection{} (Reid) Let $\Gamma$ be a Kleinian group of finite covolume.  Does $\mathcal{C}(\Gamma)$ determine $\Gamma$ amongst Kleinian groups / fundamental groups of compact 3-manifolds?

\subsection{} (Gabai, Trnkova) Let $M$ be a hyperbolic 3-manifold, $n$ - any positive integer. Does $M$ admit
an ideal triangulation with $m$ positively oriented tetrahedra, where $m\ge n$?

\subsection{} (Walsh) If $G$ is a Gromov hyperbolic group whose boundary is
homeomorphic to the limit set of a convex cocompact Kleinian group, is
$G$ virtually a convex cocompact Kleinian group? 

Known results:  When  $\partial G = S^1$ the answer is yes, and this is due to Tukia,
Gabai and Casson--Jungreis.  When $\partial G$ contains no Sierpinski carpet, the answer is again yes, and due to Haissinsky.

If $\partial G$ does contain a Sierpinski carpet, this is a conjecture of Kapovich--Kleiner, generalizing the Cannon Conjecture, which is the case $\partial G = S^2$.

\subsection{} (Walsh) What subsets of $S^2$ may occur as limit sets of convex
cocompact Kleinian groups?

\subsection{} (Walsh) For which Kleinian groups does the limit set contain
a Sierpinski carpet?  For which Kleinian groups does the limit set contain
a continuum?

\subsection{} (Walsh) Characterize limit sets of graph
Kleinian groups and iterated graph-Kleinian groups. (A \emph{graph Kleinian group} is a convex co-compact Kleinian group where the
double of the convex core is a graph manifold.  An \emph{iterated}
graph Kleinian group is one whose Bowditch decomposition (see Bowditch, \cite{Bowditch98}) contains only hanging Fuchsian and graph Kleinian pieces.

\section{Random Questions}
%%%%%%%%%%%%%%%%%%%%%%%
\subsection{} (I. Kapovitch) 
  Show that a random walk on the mapping class group gives rise to a pseudo-Anosov element whose invariant foliations have generic trivalent singularities with probability that tends to $1$ as the length of the walk tends to infinity.

\subsection{} (Maher) Start with $n$ tetrahedra, and glue faces together at random.  Note that the links of vertices in this case need not be spheres, but will be (essentially) random triangulated surfaces.  Call the resulting space a pseudo-manifold.   Investigate properties of pseudo-manifolds.

\subsection{} (Maher) Rivin \cite{Rivin14} obtains many experimental results which seem extremely regular.  This may mean there is some sort of additional structure.  Investigate.

\subsection{} (Maher) Consider the orbit of a point $x$ in Teichm\"uller space under the action of the mapping class group.

Show that the proportion of orbit points in the ball of radius $r$ in the Teichm\"uller metric which are pseudo-Anosov elements whose invariant foliation have generic trivalent singularities tends to 1 as $r \to \infty$.

Show that the proportion of orbit points in the ball of radius $r$ in the
Teichm\"uller metric which give rise to hyperbolic manifolds when used as Heegaard splitting gluing maps tends to 1 as $r \to \infty$.

\section{Geometry of Numbers / Number theory and $3$--manifolds}

\subsection{} (Long) Let $f(x) \in \mathbb{Z}[x]$ be irreducible over $\mathbb{Q}$ such that all roots are real, and suppose $f(\alpha) =0$.  Define $k = \mathbb{Q}(\alpha)$, and let $\mathcal{O}_k$ be the ring of integers. How often is $\mathcal{O}_k$ a principal ideal domain?  Infinitely often?

How often is $\mathcal{O}_k$ a Euclidean domain, with the standard norm?  For quadratic fields this happens only finitely many times.

These questions motivate the consideration of the graphs in the next question.

\newcommand{\Gprime}{\Gamma_{\mbox{\footnotesize{prime}}}}
\newcommand{\Gunit}{\Gamma_{\mbox{\footnotesize{unit}}}}
\newcommand{\Ok}{\mathcal{O}_k}
\subsection{} (Long) For $\mathcal{O}_k$ as above we can form a graph $\Gunit$ whose vertices are elements of $\Ok$, so that $\alpha$ is connected to $\beta$ by an edge if $\alpha-\beta$ is a unit.  The graph $\Gunit$ is contained in a graph $\Gprime$ with the same vertex set, but also containing edges between elements which differ by a prime.

By \cite[3.1]{LongThistlethwaite},
\[ |\mbox{clique}(\Gamma_{\mbox{\footnotesize{prime}}})| \le |\mbox{clique}(\Gamma_{\mbox{\footnotesize{unit}}})| \cdot \min\{\left|\Ok/I\right|\mid I\mbox{ prime}\}
\]
Is this an equality?

\subsection{} (Manning) If $k$ is a number field which is not totally real, is there a hyperbolic $3$--manifold with trace field $k$?  This is an old question of Neumann--Reid.  See \cite{MacReid,Neumann11} for further information.

\subsection{}\label{agoldegenerate} (Agol) Can there be a degenerate Kleinian group which is not the fiber of a fibration and
has algebraic trace field \cite{Agol15}?

\subsection{} (Schleimer) More specifically, is there a singly degenerate Kleinian group where all of the matrix entries of the group elements lie in a fixed number field? 

% One suspects that the answer is ``no''. So we relax the question and ask:

% \subsection{} (Schleimer) Is there a singly degenerate Kleinian group where the matrix entries of the generators can be obtained via a quadratically converging algorithm?

%  Note that a conjecture of McMullen (regarding a self-similarity of $\mathcal{QF}(S)$ at a point of the Bers' boundary) \cite{McMullen90} implies there is a singly degenerate Kleinian group which has a linearly converging algorithm for its generators.  See also the discussion in Section 3.7 of the book \cite{McMullen96}.

\subsection{} (McMullen) Let $M$ be a finite volume hyperbolic $3$--manifold.  Suppose $M$ contains infinitely many immersed totally geodesic surfaces.  Then is $M$ arithmetic?

%%%%%%%%end hyperbolic%%%%%%%%%%%

\section{Heegard Splittings, Multisections and Knots}
%%%%%%%%%%%%%%%%%%%%%%%
\subsection{} (Agol) Do Haken hyperbolic 3-manifolds have strongly irreducible Heegaard 
splittings? 

\subsection{} (Dunfield) If you have a hyperbolic 3-manifold with torus boundary, does its profinite completion determine whether or not it's a knot complement?

\subsection{} (Futer, Schleimer) Can the unknot be detected in polynomial time? Similarly for $S^3$.

\subsection{} (Schleimer) Is there an algorithm to detect whether a Heegaard splitting is reducible (and if so, find a reducing curve)?

\subsection{} (Schleimer) Is there a classification of strongly irreducible Heegaard splittings of a given $3$--manifold?  Update: See \url{http://arxiv.org/abs/1509.05945} for recent progress.

\subsection{} (Schleimer) Is there a fast algorithm to compute distances in the curve complex? How about distances between quasiconvex sets, such as handlebody sets?

\subsection{}\label{tillman1} (Tillmann)  Theorem:  (Gay, Kirby)
  \begin{itemize}
  \item Every closed, orientable smooth 4-manifold has a trisection.
  \item Any two trisections have a common stabilizations.
  \end{itemize}
Do higher dimensional smooth manifolds always have multisections?  Do 
What is the right generalization of ``uniqueness up to stabilization'' for multisections of smooth $n$--manifolds where $n\geq 5$?

\subsection{} (Taylor) Given a 2-component link $L$ in the 3-sphere, we can attach a band joining the components in such a way so as to end up with a knot $K$ (which depends on both the original link and the choice of band). Say that the band is ``complicated" if one of the following holds:

\begin{itemize}
\item The link $L$ is a split link and the core of the band cannot be isotoped to intersect a splitting sphere in fewer than 3 points
\item The core of the band cannot be isotoped to be disjoint from any minimal genus Seifert surface for the link.
\end{itemize}

In a recent paper I showed that a knot $K$ arising by attaching a complicated band satisfies the Cabling Conjecture. What is an example of a hyperbolic knot which cannot be created by attaching a complicated band to a 2-component link?

\subsection{} (Tillmann) Can every 2-link in $\mathbb{S}^4$ be put in bridge position with respect to the standard trisection of $\mathbb{S}^4$? (Update: The paper \url{http://arxiv.org/abs/1507.08370} gives an answer to this question.)

\end{document}